# Numerical Modeling for Simulation of Contaminant Migration of Leachate in Soil Media

**Farhikhteh Samadi[1]**

PhD Student,

Civil and Environmental Engineering Department, University of Texas at Arlington,

Arlington, TX, United States

**S. Missagh Shamshiri G**

PhD Student

Civil and Environmental Engineering Department, University of Texas at Arlington,

Arlington, TX, United States

---

[1] Corresponding Author.

E-mail address: farhikhteh.samadi@mavs.uta.edu

**Abstract**

Solid waste disposal in a landfill causes leachate production whose discharge to soil layers carries soil and leachate contaminants to the groundwater. A landfill with a capacity of 3000 ton/day and an area of 0.5 square kilometers is studied to assess and control the environmental impacts of leachate discharge on the inner layers of soil. An analytical calculation is performed, a numerical model under two different boundary conditions is developed, and their outcomes are compared to the experimental values. The landfill data is used for programming and a finite element modelling is applied. The variables considered in the model are time steps and diffusion coefficient. Results obtained from numerical computations were compared with the experimental values. The comparison demonstrates a good agreement between experimental and numerical calculations. This model can be used at fields for the landfills to simulate transport of different types of ions available in the leachates in any scale.

**Introduction**

Municipal solid waste (MSW) is continuously produced by human being all over the world. In most countries MSW is disposed of through landfills. However, leachate migration from landfills to depths of the soils causes environmental degradation and soil and groundwater contamination (Mahini and Gholamalifard, 2006). The ions existing in the leachates pass through layers of soil by different mechanisms and reach layers of soil that contain water, and eventually cause groundwater contamination (Gao et al., 2013). Pollution spreads in environment as a result of direct groundwater usage or its uptake by plants and then consumption of plants by animals and human beings (Mirbagheri and Samadi, 2013). The contaminant transport mechanisms consist of two categories: Physical migration such as advection, dispersion, diffusion and chemical transport such as sorption and precipitation (Gupta and Singh, 1997, Sharma and Lewis, 1994).

In diffusion the moving force of ions is the concentration gradient (Shackelford and Daniel, 1991). In advection process chemical ions in the soil move with flow caused by the hydraulic gradient. If the static level of leachate is higher than groundwater level, the water movement would be downward, and advection and diffusion will be in the same direction (Rowe et al., 1988). Soils that are in the path of contaminants, affect the amount and direction of the pollutants. Clay, for example, because of having an electric charge, absorbs more chemical compounds with electrical charges, than other soil types. Therefore, it is important to choose the landfill liner and layers of soils in a way that prevent the groundwater from being contaminated (Patil and Chore, 2014, Rowe and Booker, 1985, Wu and Jeng, 2017). To mitigate leachate transport from landfills to groundwater, many studies has been done. For this purpose, (Benson et al., 1999) studied different types of clay liners, (Tang et al., 2018) investigated the impact of biologically clogged barrier in landfill liners, (Reddy and Kulkarni, 2010) tried vertical wells in a bioreactor landfill and other types of barriers such as geomembrane and compacted clay layer or geomembrane and geosynthetic clay have been used (Hoor and Rowe, 2013, Varank et al., 2011).

Tracking fate of contaminants is a hard, long-term, and sometimes an impossible process. Therefore, many numerical models have been offered to be used in solving advection-dispersion reaction equations (ADRE) and simulating contaminant transport in the soil (Choi et al., 2005, Rowe et al., 1988, Shackelford and Daniel, 1991). Finite difference (FD) and finite volume (FV) methods are more popular because they are straightforward and comparatively simpler to solve compared to the existing numerical methods (Abriola, 1987, Ataie-Ashtiani et al., 1996, Mirbagheri et al., 2009, Moldrup et al., 1996).

In this paper, the goal is to solve the pollution transfer equation in the soil for a concurrent diffusion and advection using a finite difference (FD) method through MATLAB programming and under different boundary conditions and then compare the results with the ones from analytical solution and experiments with the aim of reaching an accurate estimate of final pollutants' concentration in soil layers.

**Mathematical Modelling**

One of the commonly used equations in modeling chemical transport through advection-diffusion process in the soil is used herein (Eq. (1)), which is developed based on Fick's law (Crooks and Quigley, 1984).

$$f = \theta C v - \theta D \frac{\partial C}{\partial Z} \tag{1}$$

In which $f$ is mass flux ($mg/m^2.s$) in z direction, $\theta$ is volumetric amount of water in a saturated soil and is presented in Table 1Table 1 (Reddi and Inyang, 2000). $D$ is dispersion coefficient ($m^2/s$), $v$ is Darcey Velocity ($m/s$), and $C$ is concentration of contaminant ($mg/L$).

Taking into consideration the mass balance, (Eq. (2)) it is demonstrated that (Rowe and Booker, 1985)

$$-\left(\frac{\partial f}{\partial z}\right) = \theta \left(\frac{\partial C}{\partial t}\right) + \rho K \left(\frac{\partial C}{\partial t}\right) \tag{2}$$

That can be written as

$$-\left(\frac{\partial f}{\partial z}\right) = \left(\frac{\partial C}{\partial t}\right)(\theta + \rho K) \tag{3}$$

In (Eq. (2)), K is distribution factor and is obtained from experiment ($m^3$/mg). The last term in this equation is geochemical reaction for ion exchange under equilibrium condition. But it is assumed that concentration of one exchange ion is low so that a linear relationship can be defined between the absorption of the pollutant and concentration in pore fluid (Rowe and Booker, 1985).

(Eq. (3)) can be presented as (Eq. (4)), for a saturated and homogenous soil which is at steady-state condition.

$$(\theta + \rho K_d)\frac{\partial C}{\partial t} = \theta D \frac{\partial^2 C}{\partial Z^2} - \theta v \frac{\partial C}{\partial Z} + \theta D \frac{\partial^2 C}{\partial X^2} - \theta v \frac{\partial C}{\partial X} \qquad (4)$$

Where, $K = K_d$ is distribution factor and is obtained from experiment ($m^3$/mg), t is time, $x$ and $z$ are distances along $x$ and $z$ axes respectively, $\rho$ is the bulk density of the dry soil $(mg/m^3)$. (Eq. (4)) is written for $x$ and $z$ direction and $\theta$ is considered equal in both $x$ and $z$ directions, and for solving the equations it was assumed that

$$D_x = D_z \qquad (4a)$$

$$v_x = v_z \qquad (4b)$$

Therefore, the 1-D ADRE equation in the direction of soil depth is presented as

$$R\left(\partial C / \partial t\right) = D\left(\partial^2 C / \partial Z^2\right) - v\left(\partial C / \partial Z\right) \qquad (5)$$

R is called retardation factor and can be derived from the following equation:

$$R = 1 + \left(\rho K_d / \theta\right) \qquad (6)$$

For some compounds $R = 1$ (Reddi and Inyang, 2000, Rowe and Booker, 1985).

**Model Verification**

The analytical solution to ADRE is available in the literature (Ogata, 1970, Rowe and Booker, 1985, Thongmoon and McKibbin, 2006). The method proposed by Ogata (Ogata, 1970) is commonly used in validating the results from FD models.

In this study the results from analytical solution are referred to as exact solution and are compared with the numerical solution values. The equation used for pollutant concentration at different times and distances from source of pollution, considering advection-diffusion mechanisms, is mentioned below.

$$C(z,t) = \frac{C_0}{2}\left[erfc\left(\frac{R_z - v_z t}{2\sqrt{RDt}}\right) + \exp\left(\frac{v_z z}{D}\right)erfc\left(\frac{R_z + v_z t}{2\sqrt{RDt}}\right)\right] + \frac{C_0}{2}\left[erfc\left(\frac{R_x - v_x t}{2\sqrt{RDt}}\right) + \exp\left(\frac{v_x x}{D}\right)erfc\left(\frac{R_x + v_x t}{2\sqrt{RDt}}\right)\right]$$

(7)

In this equation, $C_0$ is the surface concentration (z=0) assumed to be constant, and $z$ is the depth in which concentration is measured and $t$ represents the time. The initial and boundary conditions used in the model are as follows:

Initial condition

$C(x,z,t) = 0$ at $t = 0, z, x \geq 0$ (7a)

Boundary conditions

$C(x,z,t) = C_{i-1} = C_{i+1}$ at $x = 0$ (7b)

$C(x,z,t) = C_{i-1} = C_{i+1}$ at $x = 9$ (7c)

$C(x,z,t) = C_0$ at $z = 0$ (7d)

$\left.\dfrac{\partial C}{\partial t}\right|_{z=11} = 0$ (7e)

It is hard to control fate of contaminants and their path especially when there is interaction with soils. Therefore, the boundary conditions are unknown. To find a boundary condition that fits the FDM under conditions of the current study, the problem is solved with two different Neumann (Arendt and Warma, 2003, Dai, 2010) and reflect boundary conditions (Trefethen, 1996) with MATLAB programming. Results from MATLAB are then compared to each other as well as those from analytical solutions.

In this problem, the dimensions were considered as $x = 8$ cm and $z = 10$ cm, therefore, a $8\times10$ matrix with 9 nodes in $x$ axis and 11 nodes in $z$ axis was formed. The discretized equation for each node $(i, j)$ is written as (Ames, 2014, Anderson, 1995)

$$(\theta + \rho K)\left(\frac{C_{i,j}^{n+1} - C_{i,j}^n}{\Delta t}\right) = \theta D\left(\frac{C_{i+1,j}^n - 2C_{i,j}^n + C_{i-1,j}^n}{\Delta x^2}\right) - \theta v\left(\frac{C_{i+1,j}^n - C_{i,j}^n}{\Delta x}\right) + \theta D\left(\frac{C_{i,j+1}^n - 2C_{i,j}^n + C_{i,j-1}^n}{\Delta z^2}\right) - \theta v\left(\frac{C_{i,j+1}^n - C_{i,j}^n}{\Delta z}\right)$$

(8)

The equation is solved for $t = 0$ and all the coefficients are calculated by knowing the initial value of $C$. Then, by moving through the next time steps, and following the same procedure, the final concentration is calculated for the period that the landfill is under study.

**Site Description**

The data provided in this paper were taken from a landfill with the capacity of 3000 ton/day and area of 0.5 square kilometer. The soil type is clay with 0.3 porosity. In the experiment only KCl and the concentration of $K^+$ and $Cl^-$ ions are measured but the same model can be used for other types of contaminants. The ambient temperature of the area varies in the range of -12 to 42ºC in one year and average annual precipitation is 325 mm. The groundwater level was 40 m below the landfill. The soil characteristics that was used to simulate chemical concentration as they pass through the soil are mentioned in Table 1. The depth of the soil used in the exact method was from $z = 0$ to $z = 11$ cm. The landfill was under study for 10 years and 10-year annual measured mean concentrations is used.

Table 1. Soil characteristic in the landfill area

| Parameters | Values |
| --- | --- |
| Initial concentration for $Cl^-$ and $K^+$ ($C_0$, mg/L) | 675 |
| Study Period (Days) | 100 |
| Soil Porosity ($\theta$) | 0.30 |
| diffusion-dispersion coefficients [a] (D, m$^2$/a) | 0.02 |
| Back ground contaminant level ($K^+$ and $Cl^-$, mg/L) | Trace |

[a] (Rowe, et al., 1988)

**Results and Discussions**

Results from the landfill are presented in Figure 1. It can be seen that as KCl travels through the soil, it is absorbed by soil and its concentration decreases.

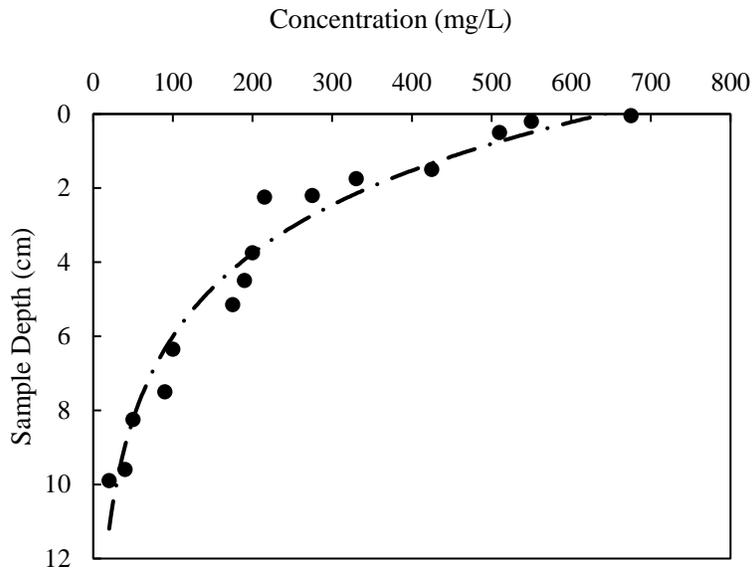

Figure 1. Concentration profile of KCl in soil depth

As mentioned before the problem was solved for Neumann boundary condition. Figure 2 illustrates the $K^+$ ion concentration profile using MATLAB programming and Neumann boundary condition.

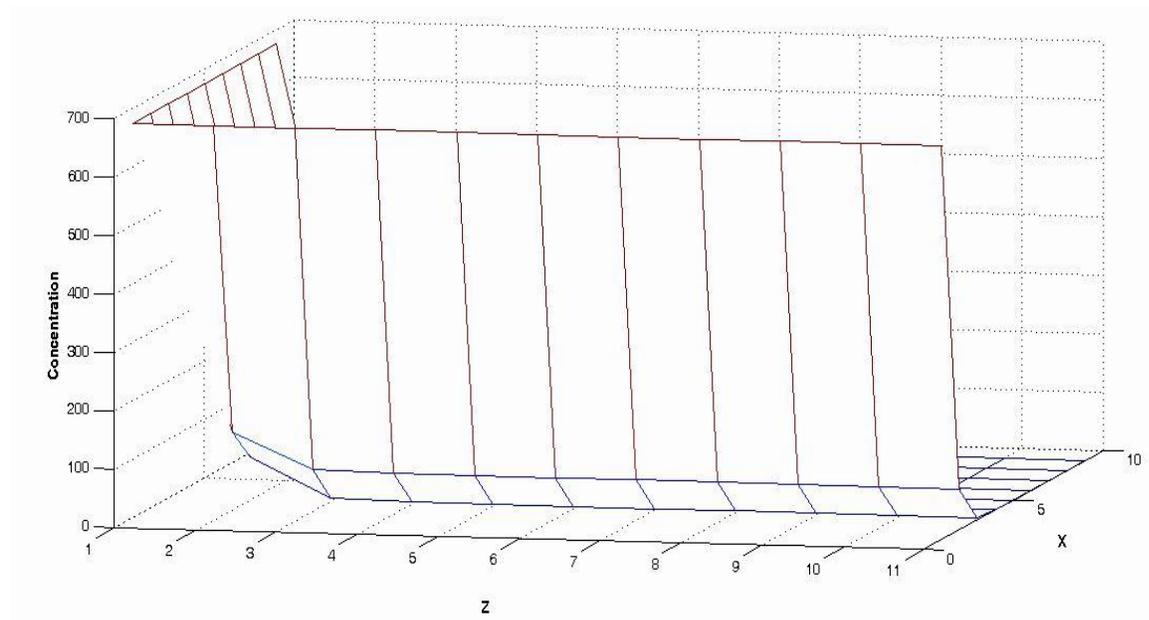

**Figure 2. Concentration profile for potassium ion using Neumann boundary condition**

In order to make a better comparison, the FD method is solved for the two boundary conditions and the results are presented in Figure 3. The results are shown for node *x = 5* and *z = 11*.

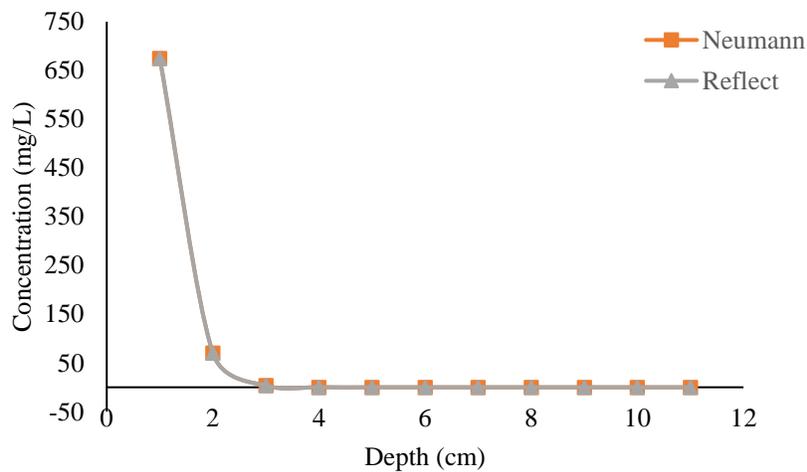

**Figure 3. Concentration profile for node x = 5, z = 11 solved with Neumann and reflect method**

It can be seen in Figure 2 and Figure 3 that concentration in depth has a reducing trend, as expected from fig1. Meaning that the soil is transferring a very small amount of pollutant to lower layers.

In reflect condition it was assumed that the soil environment has the same conditions in a wide range.

**Results of FD solution for K⁺ ion was solved through exact model, Neumann and reflect boundary condition for different days. (a) t = 1 day            (b) t = 50 day**

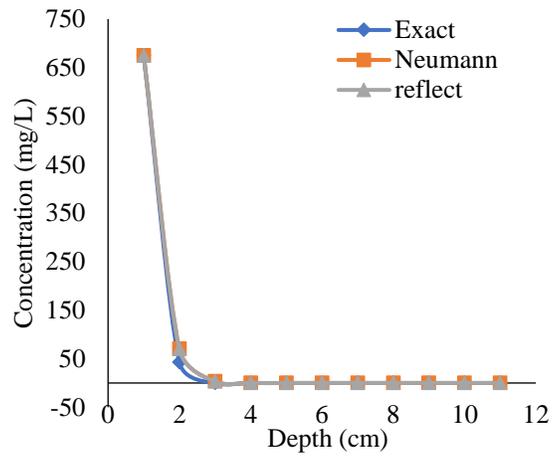

**(c) t = 100 day**

Figure 4 a, b and c show results for day 1, day 50 and day 100 of the experiment, respectively.

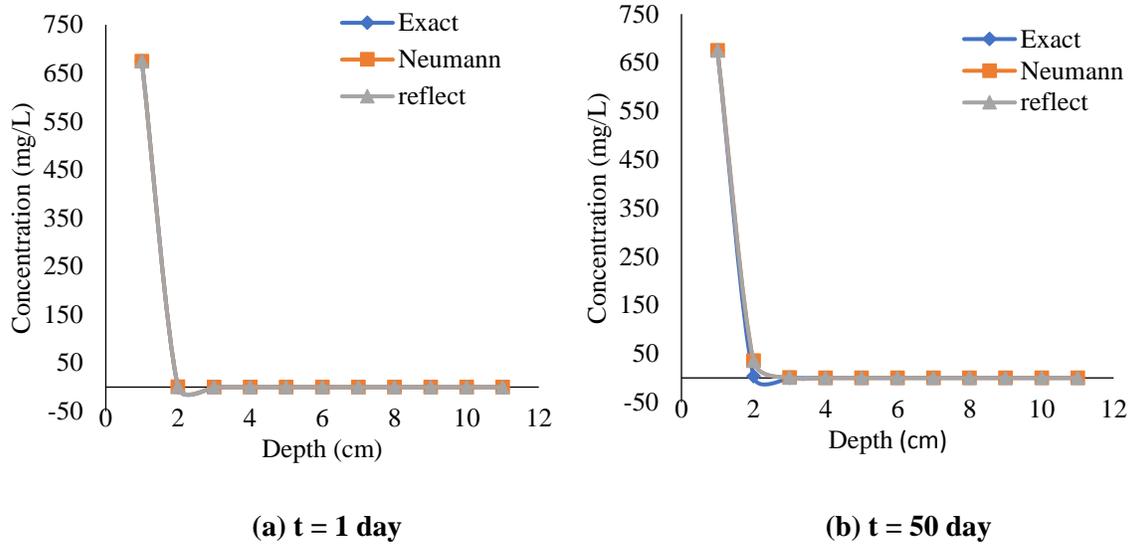

(a) t = 1 day

(b) t = 50 day

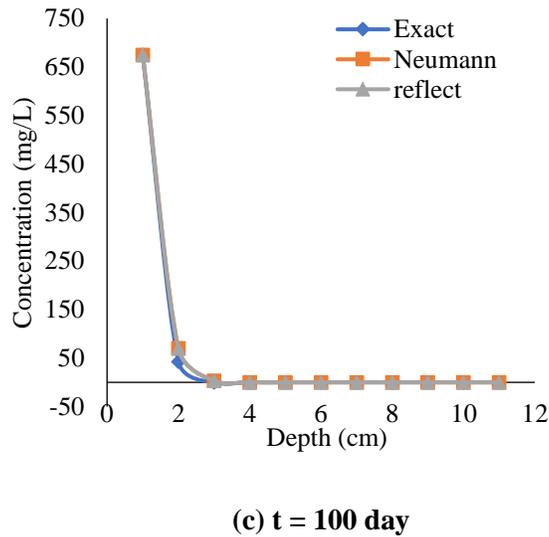

(c) t = 100 day

**Figure 4. FD solution for K$^+$ concentration-depth profile with Neumann and reflect boundary condition, and exact solution for (a) t=1, (b) t=50, and (c) t= 100 (day)**

Figure 4 illustrates that there is sudden change in concentration of cation from depth 0 to 2 cm and then a moderate reduction of concentration is observed from $z = 2$ to $z = 10$. Therefore, it can be concluded that according to FD model most of the cation is absorbed in the depths closer to the soil surface. Also, in Figure 4 (a) and (b) there is a slight difference between exact solution and FD models but at $t = 100$ day, they all follow the same trend and the graphs match each other.

The FD model using Neumann boundary condition is then solved for Cl$^-$ ion via MATLAB program. The graph is provided in Figure 5.

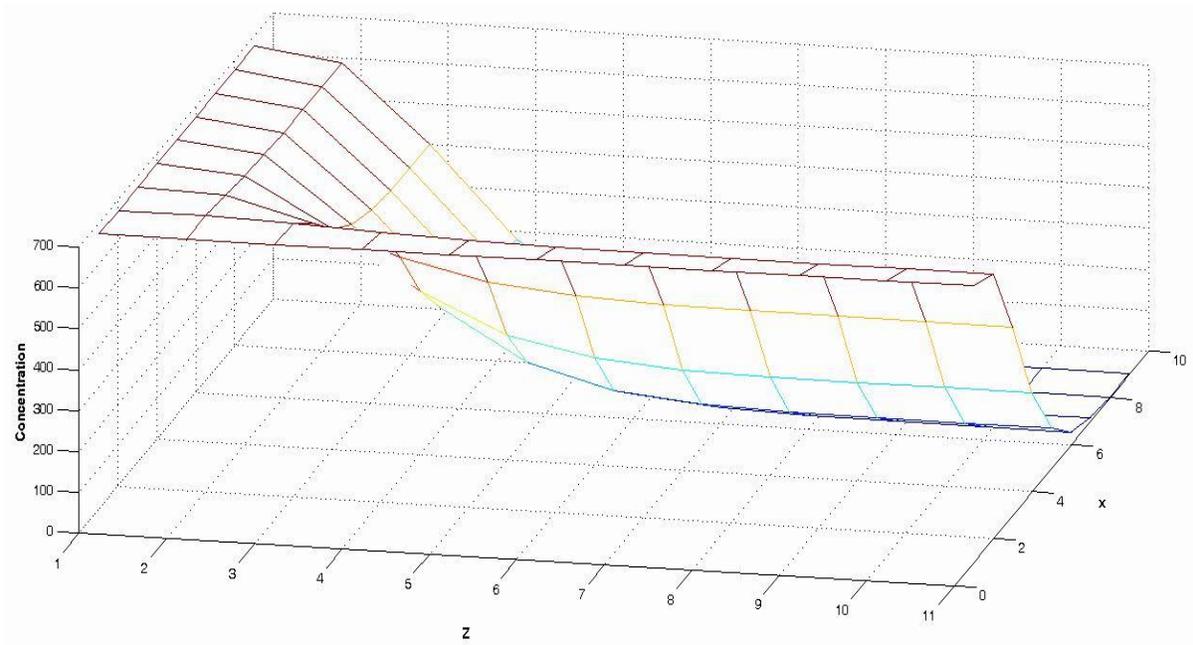

**Figure 5. Concentration profile for chloride ion using Neumann boundary condition**

The FD simulation results of the model under Neumann and reflect boundary conditions are shown in Figure 6. The anion behavior is totally different from cation behavior, in this figure. The concentration of chloride ion does not change up to $z = 3$ cm while it was reduced to a great extent for potassium ion in the same depth. After $z = 3$ cm, Cl$^-$ reduces gradually, and in $z = 10$ it reaches the trace amount, under Neumann boundary condition. However, according to reflect boundary condition for the chloride ion concentration, soil is not capable of absorbing much contaminant and the final concentration of Cl$^-$ gets only as low as 150 mg/L.

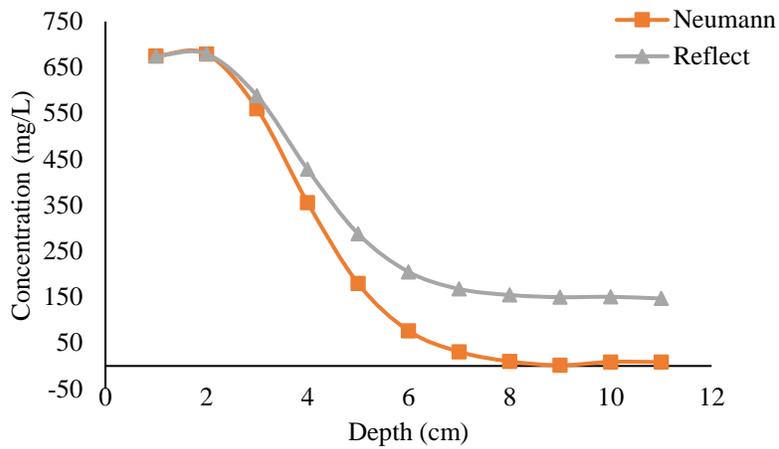

**Figure 6. Concentration profile comparison for node x = 5, z = 11 solved with Neumann and reflect method**

(a) In the next figure, the same as Figure 4, a comparison was made between analytical solution and FD models for chloride ion at different time steps of *t* = 1 day, *t* = 50 day, and *t* = 100 day and the results are provided in  t = 1 day
(b) t = 50 day

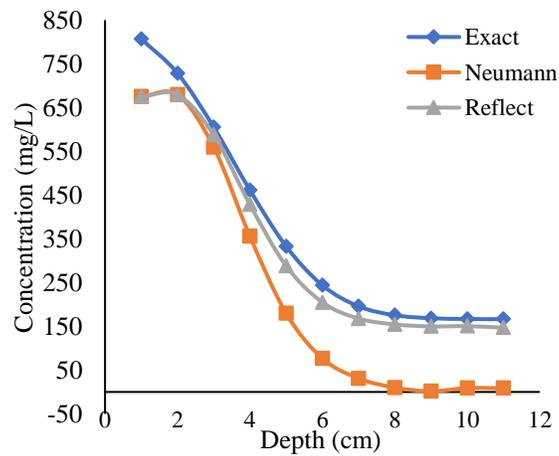

(c) t = 100 day

Figure 7 (a), (b) and (c), respectively.

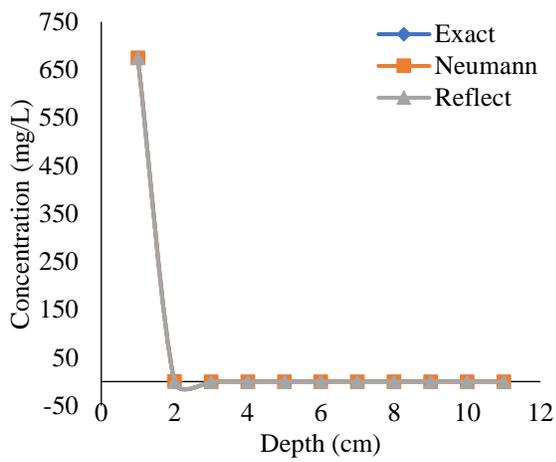
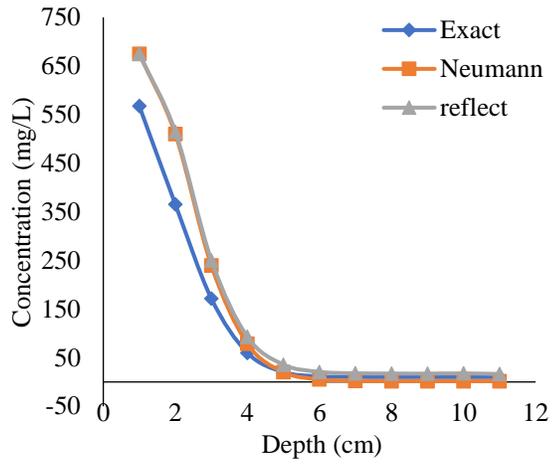

(b) t = 1 day

(b) t = 50 day

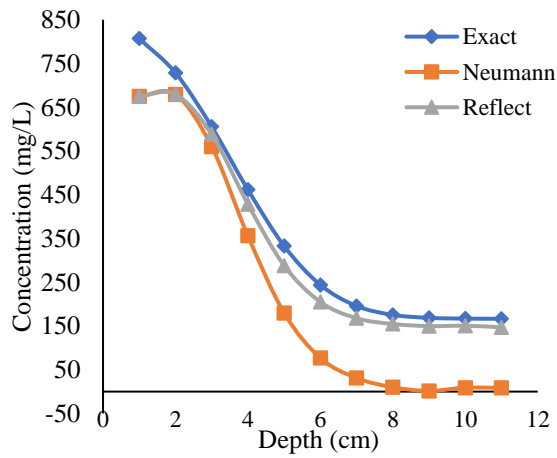

(c) t = 100 day

**Figure 7. FD solution for Cl⁻ concentration-depth profile with Neumann and reflect boundary condition and exact solution for (a) t = 1, (b) t = 50, and (c) t = 100 (day)**

In Figure 7 (c), reflect boundary condition gave closer values to analytical solution, than the Neumann boundary condition.

Moreover, in Figure 7, using larger time steps for the FD model results in a bigger difference between the exact and FD solutions. The exact method is calculated through eq. 4, but for numerical solution a mesh system is used. When either the meshes or the time steps get smaller, the numerical results will get closer to the analytical ones, and this can go on up to a point that mesh sizes do not make any difference in the results. This is called grid independency. The problem is solved numerically for different time steps changing from t = 100 day to t = 0.01 day, and the results of concentration in depth are presented in

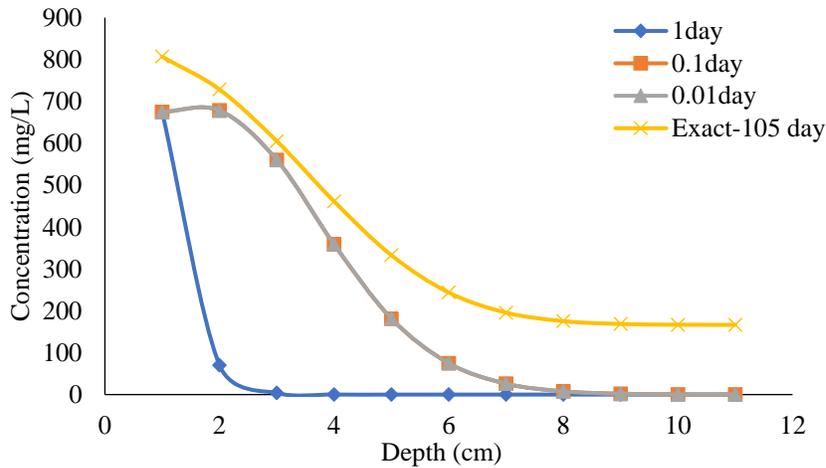

**Figure 8**. It shows a significant difference between the changes in the concentration in 100 day and 0.1 day but the graphs for t = 0.1 day and t = 0.01 day do not show a notable difference.

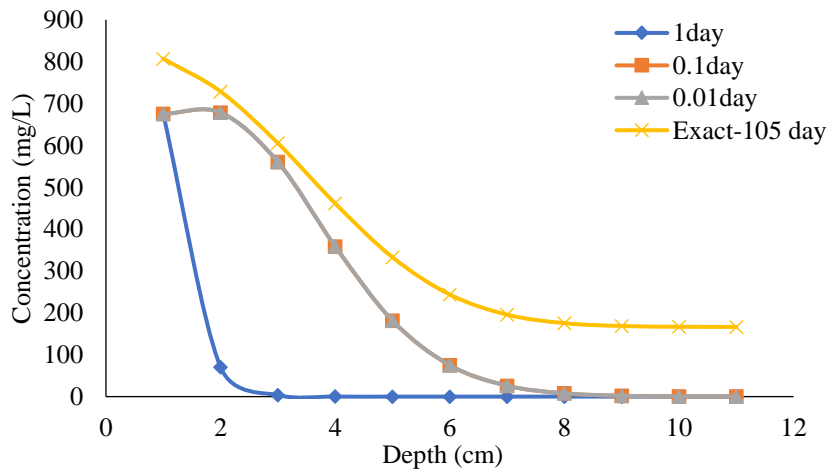

**Figure 8. Comparison of Exact and FD results at different time steps with Neumann boundary condition for Cl- ion**

Furthermore, the mesh sizes changed in both x and z direction and the problem was solved to compare numerical and analytical solutions. The results are provided in Figure 9. It was noticed that the results are more scattered if the mesh was too coarse or too fine; therefore, the element size of x, z = 1 cm is chosen for the numerical solution.

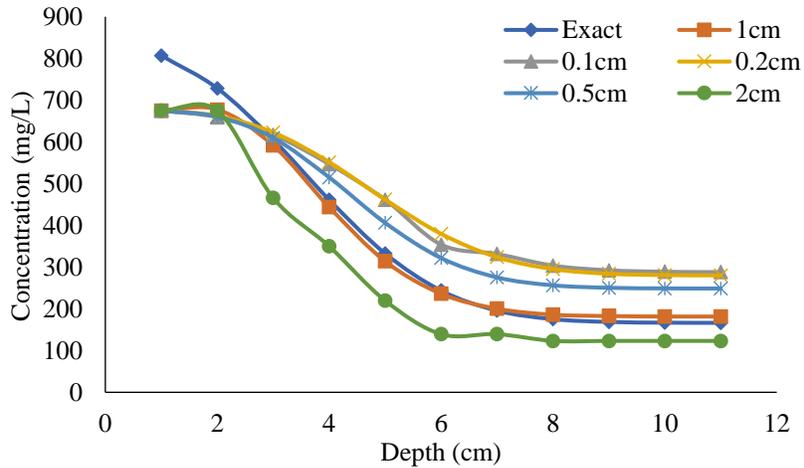

**Figure 9. Comparison of Exact and FD solution at different mesh sizes with Neumann boundary condition for Cl⁻ ion**

A diffusion coefficient of 0.02 m$^2$/a was used in the previous studies available in the literature, while the minimum value for this parameter was 0.018 m$^2$/a, suggested by (Rowe and Booker, 1985). Hence, in the current study, two different FD models are solved for both $D = 0.02$ m$^2$/a, and $D = 0.018$ m$^2$/a. The results are presented in Figure 10, showing that increasing the diffusion coefficient from 0.018 m$^2$/a to 0.02 m$^2$/a makes a slight change in the values, but the overall trend remains the same.

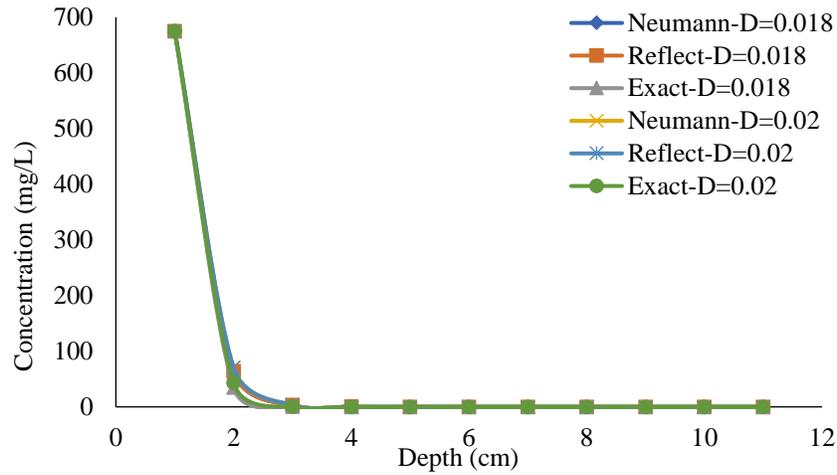

**Figure 10. Comparison of the effect of change in diffusion coefficient in concentration profile for Cl⁻ ion**

**Conclusion**

A finite difference model was developed to predict the concentration profile of contaminants in a clayey soil under the landfill. A landfill was studied, and the gained data were used for numerical modelling. Two types of boundary conditions of Neumann and reflect were applied to the problem domain. Analytical solution was also performed, and the results were used to validate the numerical solution's results. It was observed that cation concentration in soil reduces significantly at layers close to the soil surface, but the anion concentration changes gradually along the soil depth. In this study, $K^+$ and $Cl^-$ ions were examined, and the model can be used for other types of contaminants. Although, the soil showed a strong capacity in absorbing pollutants and preventing them from contaminating the groundwater, it is recommended to utilize the model in field scale problems.

**Notation List**

| | |
|---|---|
| C | concentration of contaminant, mg/m³ |
| D | dispersion coefficient, m²/s |
| f | mass flux, mg/m$^2$s |
| K | distribution factor, m$^3$/kg |
| R | retardation factor, dimensionless |
| t | time, s |
| v | Darcy velocity, m/s |
| X | distance in Cartesian coordinate, m |
| Z | depth in Cartesian coordinate, m |

**Greek Letters and Subscripts/Superscript**

| | |
|---|---|
| θ | soil porosity, dimensionless |
| i | x direction elements counter |
| j | z direction elements counter |
| n | temporal element number |
| x | in x direction |
| z | in z direction |